# Nikol'skii–Type Inequalities for Entire Functions of Exponential Type in Lorentz–Zygmund Spaces

**Leo R. Ya. Doktorski**[1]

**Abstract**

Nikol'skii–type inequalities for entire functions of exponential type on $\mathbb{R}^n$ for the Lorentz–Zygmund spaces are obtained. Some new limiting cases are examined. Application to Besov–type spaces of logarithmic smoothness is given.

## 1. Introduction

Let a function $f \in S'(\mathbb{R}^n)$ be such that $supp\hat{f}$ is compact. Due to Paley–Wiener–Schwartz theorem (see [23] and [35, 2.2.1]), it is an entire function of exponential type on $\mathbb{C}^n$. Let additionally $f \in L_q(\mathbb{R}^n)$. For such functions Nikol'skii has proved [24, 25] that for some $\omega$ (depending on $supp\hat{f}$)

$$\|f\|_p \leq C\omega^{n\left(\frac{1}{q}-\frac{1}{p}\right)}\|f\|_q, \tag{1}$$

where $1 \leq q < p \leq \infty$, $\|*\|_q$ is the usual norm on the Lebesgue spaces $L_q$ and $C$ does not depend on $f$. The proof in [24, 25] is based on Bernstein's inequality.

Inequalities between different (quasi-)norms of the same function are known as Nikol'skii–type inequalities. They play a crucial role in many areas of mathematics, e.g. theory of approximation, theory of functions of several variables, and functional analysis (embedding theorems for Besov spaces). Nessel and Wilmes [23] extended (1) for $0 < q < p \leq \infty$. Ditzian and Prymak [6] used Bernstein–type inequalities and generalised (1) to the Lorentz spaces. Nikol'skii–type inequality for entire functions from the Lorentz–Zygmund spaces is presented and used in [18]. The classical Nikol'skii inequality for the trigonometric polynomials $T_n$ on $[0,1]$ of degree at most *n* can be written as [24, 25]

$$\|T_n\|_p \leq Cn^{\frac{1}{q}-\frac{1}{p}}\|T_n\|_q. \tag{2}$$

This variant of Nikol'skii's inequality is also intensively studied in [6, 18, 23, 32]. In [11, 12], (2) is generalised to the Lorentz–Zygmund spaces; some new limiting cases are obtained. In numerous other articles different sets of functions, domains and measures are explored. For further information about Nikol'skii–type inequalities and theirs applications we refer to [4, 6, 7, 8, 11, 12, 13, 14, 18, 23, 28, 32, 33, 34] and references therein.

The principal aim of this paper is to extend the inequality (1) to the Lorentz–Zygmund spaces. The technique we apply relies on properties of the Lorentz–Zygmund spaces and the ones of the Fourier transform. Additionally, we use Hardy-type inequality, some facts from the

---

[1] Fraunhofer Institute of Optronics, System Technologies and Image Exploitation IOSB, Department Object Recognition, Gutleuthausstr. 1, 76275 Ettlingen, Germany. ORCID iD: 0000-0001-5977-1000

Correspondence should be addressed to Leo Doktorski; doktorskileo@gmail.com or leo.doktorski@iosb.fraunhofer.de

real interpolation method, and the observation that the power of an entire function is also an entire function [23, 34].

This paper is organized as follows. Section 2 contains necessary notations and definitions. Main results of this contribution are formulated and proved in Section 3. In Section 4, we define some Besov–type spaces of logarithmic smoothness (as well homogeneous counterparts) and apply main results to establish embeddings between these spaces.

## 2. Preliminaries

Throughout the paper we write $X \subset Y$ for two quasi-normed spaces $X$ and $Y$ to indicate that $X$ is continuously embedded in $Y$. The notation $X = Y$ means that $X \subset Y$ and $Y \subset X$. In this case $X$ and $Y$ are equal as sets and as linear spaces and their (quasi-)norms are equivalent. We adopt the convention that $\frac{1}{\infty} = 0$ and define $p'$ by $\frac{1}{p'} + \frac{1}{p} = 1$ ($1 \leq p \leq \infty$). $\|*\|_{q,(a,b)}$ is the usual (quasi-)norm on the Lebesgue space $L_q$ on the interval $(a, b)$ ($0 < q \leq \infty, 0 \leq a < b \leq \infty$). We shall use bold capital Greek letters for ordered pairs of real numbers: anywhere below $\mathbf{A} = (\alpha_0, \alpha_\infty), \mathbf{B} = (\beta_0, \beta_\infty), \mathbf{\Gamma} = (\gamma_0, \gamma_\infty), \mathbf{\Delta} = (\delta_0, \delta_\infty) \in \mathbb{R}^2$. For those pairs we use usual vector operations and the conventions $\mathbf{A} + \sigma = (\alpha_0 + \sigma, \alpha_\infty + \sigma)$ and put $\overleftarrow{\mathbf{A}} = (\alpha_\infty, \alpha_0)$. We also write $\mathbf{A} = 0$, when $\alpha_0 = \alpha_\infty = 0$, and $\mathbf{A} < \mathbf{B}$, when $\alpha_0 < \beta_0$ and $\alpha_\infty < \beta_\infty$. We use the abbreviation $l(t) = 1 + |\ln t|$ and $ll(t) = l(l(t)), 0 < t < \infty$. For an ordered pair $\mathbf{A}$ a broken-logarithmic function is defined by
$$l^{\mathbf{A}}(t) = \begin{cases} (1 + |\ln t|)^{\alpha_0}, & 0 < t \leq 1, \\ (1 + |\ln t|)^{\alpha_\infty}, & t > 1 \end{cases}$$
and analogously for $ll^{\mathbf{A}}(t)$.

As usual, by $S \equiv S(\mathbb{R}^n)$ and by $S' \equiv S'(\mathbb{R}^n)$ we denote the space of rapidly decreasing functions and the space of tempered distributions on $\mathbb{R}^n$, correspondently.

### 2.1. Real interpolation method

Let $\bar{X} := (X_0, X_1)$ be a compatible couple of (quasi-) Banach spaces such that $X_0 \cap X_1 \neq \{0\}$. The Peetre's $K$-functional on $X_0 + X_1$ is given by
$$K(t, f; \bar{X}) \equiv K(t, f) := \inf(\|f_0\|_{X_0} + t\|f_1\|_{X_1}) \quad (f = f_0 + f_1, f_i \in X_i, i=0,1, t > 0).$$
The classical (Lions-Peetre) scale of interpolation spaces $\bar{X}_{\theta,q} \equiv (X_0, X_1)_{\theta,q}$ ($0 < \theta < 1$) is defined via the (quasi-)norms
$$\|f\|_{\theta,q} := \|u^{-\theta-1/q} K(u, f)\|_{q,(0,\infty)}. \tag{3}$$
We will often use the following properties of the scale $\bar{X}_{\theta,q}$.

**Lemma 1** Let $0 < \theta < 1$ and $0 < q \leq \infty$.
  (i) $(*, *)_{\theta,q}$ is an exact interpolation functor of exponent $\theta$.
  (ii) There exists a positive number $C_{\theta,q}$ such that for all $f \in X_0 \cap X_1$
$$\|f\|_{\theta,q} \leq C_{\theta,q} \|f\|_{X_0}^{1-\theta} \|f\|_{X_1}^{\theta}.$$

For further information about basic properties of the $K$-functional and the real interpolation method we refer to e.g. [2, 3, 35]. The definition of the classical functor $(*, *)_{\theta,q}$ requires $0 < \theta < 1$ or $\theta \in \{0,1\}$ for $q = \infty$, yet extreme cases $\theta \in \{0,1\}$ and $0 < q < \infty$ play an important role in certain questions in analysis. One of the possibilities to consider these extreme cases is to involve an additional logarithmic factor $l^{\mathbf{A}}(t)$ in the formula (3).

**Definition 2** (See, e.g., [9, 16, 17, 19].) Let $0 \leq \theta \leq 1$, $0 < q \leq \infty$, and $\mathbf{A} = (\alpha_0, \alpha_\infty)$. We put

$$\bar{X}_{\theta,q;\mathbf{A}} \equiv (X_0, X_1)_{\theta,q;\mathbf{A}} := \left\{ f \in X_0 + X_1 : \|f\|_{\theta,q;b} = \left\| u^{-\theta - 1/q} l^{\mathbf{A}}(t) K(u, f) \right\|_{q,(0,\infty)} < \infty \right\}.$$

Note that $(*,*)_{\theta,q;0} = (*,*)_{\theta,q}$. We will use the following properties of the scale $(X_0, X_1)_{\theta,q;\mathbf{A}}$. See e.g. [17, 19].

**Lemma 3** Let $0 \leq \theta \leq 1$, $0 < q \leq \infty$, and $\mathbf{A} = (\alpha_0, \alpha_\infty)$.

(i) $X_0 \cap X_1 \subset \bar{X}_{\theta,q;\mathbf{A}} \subset X_0 + X_1$ if and only if one of the following conditions holds true:

$$\begin{cases} 0 < \theta < 1, \\ \theta = 0, \alpha_\infty + \frac{1}{b} < 0, \\ \theta = 0, b = \infty, \alpha_\infty = 0, \\ \theta = 1, \alpha_0 + \frac{1}{b} < 0, \\ \theta = 1, b = \infty, \alpha_0 = 0. \end{cases} \quad (4)$$

(ii) If none of the conditions (4) holds then $\bar{X}_{\theta,q;\mathbf{A}} = \{0\}$.
(iii) If one of the conditions (4) holds then $(*,*)_{\theta,q;\mathbf{A}}$ is an exact interpolation functor.
(iv) $(X_1, X_0)_{\theta,q;\mathbf{A}} = (X_0, X_1)_{1-\theta,q;\bar{\mathbf{A}}}$.

### 2.2. Slowly varying functions

In this section, we summarise some properties of slowly varying functions, which will be required later. For more details, we refer to e.g. [19].

**Definition 4** We say that a positive Lebesgue measurable function $b$ is slowly varying on $(0,\infty)$, notation $b \in SV$, if for each $\varepsilon > 0$ the function $t^\varepsilon b(t)$ is equivalent to an increasing function while the function $t^{-\varepsilon} b(t)$ is equivalent to a decreasing function.

The functions of the form $l^{\mathbf{A}}(t)$ are examples of slowly varying functions. We need following properties of the $SV$-functions.

**Lemma 5** Let $b \in SV$, $\sigma > 0$, $0 < r \leq \infty$, $\lambda \in (-\infty, \infty)$, and $t \in (0, \infty)$. Then

(i) $ct^\sigma b(t) \leq \left\| s^{\sigma - \frac{1}{r}} b(s) \right\|_{r,(0,t)} \leq Ct^\sigma b(t)$.

(ii) $ct^\lambda b(t) \leq \left\| u^{\lambda - \frac{1}{r}} b(u) \right\|_{r,\left(\frac{t}{2},t\right)} \leq Ct^\lambda b(t)$.

The constants $c$ and $C$ in both inequalities do not depend on $t$.

### 2.3. Lorentz–Zygmund spaces

We consider (equivalence classes of) complex-valued measurable functions $f$ on $\mathbb{R}^n$ ($n=1,2,\ldots$) with the standard Lebesgue measure $\mu$. As conventional [1, 2, 27, 31], $f^*(t)$ ($t>0$) denotes the non-increasing rearrangements of a function $f$. The maximal function of $f^*$ is defined by $f^{**}(t) := \frac{1}{t} \int_0^t f^*(u) du$.

**Definition 6** Let $M$ be a measurable set in $\mathbb{R}^n$ and $0<p,b\leq\infty$. The Lorentz–Zygmund space $L_{p,b;\mathbf{A}}(M)$ is defined by

$$L_{p,b;\mathbf{A}}(M) := \left\{ f \in \mathfrak{M}(M) : \|f\|_{p,b;\mathbf{A}} := \left\| t^{\frac{1}{p} - \frac{1}{b}} l^{\mathbf{A}}(t) f^*(t) \right\|_{b,(0,\mu(M))} < \infty \right\},$$

where $\mathfrak{M}(M)$ is the set of all $\mu$-measurable functions on $M$.

It is known [27] that $L_{p,b;\mathbf{A}}(M)$ is not trivial, that is not equal to $\{0\}$, iff one of the following conditions holds:

$$\begin{cases} p < \infty; \\ p = \infty, \ 0 < b < \infty, \alpha_0 + \frac{1}{b} < 0; \\ p = \infty, \ 0 < b \leq \infty, \alpha_0 \leq 0. \end{cases} \quad (5)$$

If $\mu(M) < \infty$, then $f^*(t) = 0$, for $t > \mu(M)$. We write $L_{p,b;\mathbf{A}}$ for $L_{p,b;\mathbf{A}}(\mathbb{R}^n)$. $L_{p,b;(0,0)}$ coincides with the Lorentz space $L_{p,b}$ and $L_{p,p;(0,0)}$ – with the Lebesgue space $L_p$. We put $\|f\|_{p,b} = \|f\|_{p,b;(0,0)}$ and $\|f\|_p = \|f\|_{p,p}$. Note that (see e.g. [16, 17, 19])

$$L_{p,b;\Gamma} = \left(L_{p_0,b_0;\mathbf{A}_0}, L_{p_1,b_1;\mathbf{A}_1}\right)_{\theta,b;\mathbf{A}}, \quad (6)$$

where $0<\theta<1$, $0 < p_0 < p_1 < \infty$, $\frac{1}{p} = \frac{(1-\theta)}{p_0} + \frac{\theta}{p_1}$, and $\Gamma = \mathbf{A} + (1-\theta)\mathbf{A}_0 + \theta\mathbf{A}_1$. Formula (6) is also true if $p_1 = b_1 = \infty$ and $\mathbf{A}_1 = 0$. Additionally, by [17, 19]

$$L_{\infty,b;\mathbf{A}} = \left(L_{p,c;\mathbf{B}}, L_\infty\right)_{1,b;\mathbf{A}}, \quad \text{provided } 0<p<\infty \text{ and } \alpha_0 + \frac{1}{b} < 0. \quad (7)$$

More information about Lorentz–Zygmund spaces can be found in [1, 17, 19, 27, 31]. The next lemma plays a crucial role in the proof of main theorems.

**Lemma 7** (Cf. [27, Theorem 4.6].) Let $M$ be a measurable set in $\mathbb{R}^n$, $0 < \mu(M) < \infty$, $0<p<q\leq\infty$, and $0<b,c\leq\infty$. Then for all $f \in L_{q,c;\mathbf{B}}(M)$

$$\|f\|_{p,b;\mathbf{A}} \leq C\mu(M)^{\left(\frac{1}{p}-\frac{1}{q}\right)} \left(l\big(\mu(M)\big)\right)^{\mathbf{A}-\mathbf{B}} \|f\|_{q,c;\mathbf{B}},$$

where $C$ does not depend on $f$ (and $M$).

*Proof.* Below in this proof, we denote by $C$ possibly different constants which depend only on $p, b, \mathbf{A}, q, \mathbf{c}$, and $\mathbf{B}$ and do not depend on $f$ and $M$. Because $\frac{1}{p} - \frac{1}{q} > 0$, Lemma 5 (i) implies that for all $t \leq \mu(M)$

$$t^{\frac{1}{p}-\frac{1}{q}} l^{\mathbf{A}-\mathbf{B}}(t) \leq C\mu(M)^{\left(\frac{1}{p}-\frac{1}{q}\right)} \left(l\big(\mu(M)\big)\right)^{\mathbf{A}-\mathbf{B}}. \quad (8)$$

**Case 1:** $q<\infty$, $b\leq c\leq\infty$. Using Hölder's inequality with exponents $\frac{c}{b}$ and $\frac{c}{c-b}$ we get

$$\|f\|_{p,b;\mathbf{A}} = \left\|t^{\frac{1}{p}-\frac{1}{b}} l^{\mathbf{A}}(t) f^*(t)\right\|_{b,(0,\mu(M))} \leq \|f\|_{q,c;\mathbf{B}} \left\|t^{\frac{1}{p}-\frac{1}{q}-\frac{(c-b)}{bc}} l^{\mathbf{A}-\mathbf{B}}(t)\right\|_{\frac{bc}{(c-b)},(0,\mu(M))}.$$

(If $b=c=\infty$, we take here $\frac{bc}{c-b} = \infty$.) Using Lemma 5 (i) for the last factor, we derive that

$$\|f\|_{p,b;\mathbf{A}} \leq C\mu(M)^{\left(\frac{1}{p}-\frac{1}{q}\right)} \left(l\big(\mu(M)\big)\right)^{\mathbf{A}-\mathbf{B}} \|f\|_{q,c;\mathbf{B}}. \quad (9)$$

**Case 2:** $q<\infty$, $c<b=\infty$. Since $f^*$ is non-increasing, using Lemma 5 (i), we conclude that

$$t^{\frac{1}{q}} l^{\mathbf{B}}(t) f^*(t) \leq C f^*(t) \left\|s^{\frac{1}{q}-\frac{1}{c}} l^{\mathbf{B}}(s)\right\|_{c,(0,t)} \leq C \left\|s^{\frac{1}{q}-\frac{1}{c}} l^{\mathbf{B}}(s) f^*(s)\right\|_{c,(0,t)} \leq C\|f\|_{q,c;\mathbf{B}}. \quad (10)$$

Hence, (8) and (10) imply

$$\|f\|_{p,\infty;\mathbf{A}} = \sup_{0<t<\mu(M)} t^{\frac{1}{p}} l^{\mathbf{A}}(t) f^*(t) \leq \sup_{0<t<\mu(M)} t^{\frac{1}{p}-\frac{1}{q}} l^{\mathbf{A}-\mathbf{B}}(t) \sup_{0<t<\mu(M)} t^{\frac{1}{q}} l^{\mathbf{B}}(t) f^*(t)$$

$$\leq C\mu(M)^{\left(\frac{1}{p}-\frac{1}{q}\right)} \left(l\big(\mu(M)\big)\right)^{\mathbf{A}-\mathbf{B}} \|f\|_{q,c;\mathbf{B}}. \quad (11)$$

**Case 3:** $q<\infty$, $c<b<\infty$. Let $p<r<q$. Using (9) and (11), we obtain

$$\|f\|_{p,b;\mathbf{A}} \leq C\mu(M)^{\left(\frac{1}{p}-\frac{1}{r}\right)} \left(l\big(\mu(M)\big)\right)^{\mathbf{A}} \|f\|_{r,\infty} \leq C\mu(M)^{\left(\frac{1}{p}-\frac{1}{q}\right)} \left(l\big(\mu(M)\big)\right)^{\mathbf{A}-\mathbf{B}} \|f\|_{q,c;\mathbf{B}}.$$

**Case 4:** $q=\infty$ and $b = c$, $\mathbf{A} = \mathbf{B}$. Due to (5), consider two subcases. If $0<c<\infty$ and $\beta_0 + \frac{1}{c} < 0$ we derive that

$$\|f\|_{p,c;\mathbf{B}} = \left\{\int_0^{\mu(M)} \left(t^{\frac{1}{p}-\frac{1}{c}} l^{\mathbf{B}}(t) f^*(t)\right)^c dt\right\}^{\frac{1}{c}}$$

$$\leq \sup_{0<t<\mu(M)} t^{\frac{1}{p}} \cdot \left\{\int_0^{\mu(M)} \left(t^{-\frac{1}{c}} l^{\mathbf{B}}(t) f^*(t)\right)^c dt\right\}^{\frac{1}{c}} = \mu(M)^{\frac{1}{p}} \|f\|_{\infty,c;\mathbf{B}}.$$

Analogously, if $c = \infty$ and $\beta_0 \leq 0$, we have

$$\|f\|_{p,\infty;\mathbf{B}} = \sup_{0<t<\mu(M)} \left(t^{\frac{1}{p}} l^{\mathbf{B}}(t) f^*(t)\right)$$

$$\leq \mu(M)^{\frac{1}{p}} \sup_{0<t<\mu(M)} \left(l^{\mathbf{B}}(t) f^*(t)\right) = \mu(M)^{\frac{1}{p}} \|f\|_{\infty,\infty;\mathbf{B}}.$$

**Case 5:** $q=\infty$ and arbitrary triple $(p, b, \mathbf{A})$. We take some $r$, such that $p<r<\infty$. Using appropriate of the previous inequalities, we get

$$\|f\|_{p,b;\mathbf{A}} \leq C\mu(M)^{\left(\frac{1}{p}-\frac{1}{r}\right)} \left(l(\mu(M))\right)^{\mathbf{A}-\mathbf{B}} \|f\|_{r,c;\mathbf{B}} \leq C\mu(M)^{\frac{1}{p}} \left(l(\mu(M))\right)^{\mathbf{A}-\mathbf{B}} \|f\|_{\infty,c;\mathbf{B}}.$$

This completes the proof. □

**Lemma 8** Let $0<p,b\leq\infty$, $f \in L_{p,b;\mathbf{A}}$, and $k = 2,3,\ldots$. Then

$$\left\{\|f\|_{p,b;\mathbf{A}}\right\}^k = \|f^k\|_{p/k,b/k;k\mathbf{A}}.$$

### 2.4. Fourier transform

Denote by $\mathcal{F}f \equiv \hat{f}$ the Fourier transform given for a function $f \in S(\mathbb{R}^n)$ by

$$\mathcal{F}f(\xi) = \int e^{-2\pi i x \xi} f(x) dx.$$

By $\mathcal{F}^{-1}$ we denote the inverse Fourier transform. Let $\mathcal{J}$ be either $\mathcal{F}$ or $\mathcal{F}^{-1}$. It is known that

$$\|\mathcal{J}: L_1 \to L_\infty\| = 1 \tag{12}$$

and

$$\|\mathcal{J}: L_2 \to L_2\| = 1.$$

The usual real interpolation argument with $0<\theta<1$ yields

$$\|\mathcal{J}: L_{p,b;\mathbf{A}} \to L_{p',b;\overline{\mathbf{A}}}\| \leq C, \qquad 1<p<2. \tag{13}$$

Using [19] and [27 Theorem 3.8 (ii)], we get the following results for limiting case $\theta=0$:

$$\|\mathcal{J}: L_{1,1;\mathbf{A}+\mathbf{1}} \to L_{\infty,1;\overline{\mathbf{A}}}\| \leq C, \qquad \text{provided that } \alpha_\infty + 1 < 0 < \alpha_0 + 1. \tag{14}$$

The constants $C$ in (13) and (14) depend only on space's parameters $p, b, \mathbf{A}$.

### 2.5. Some auxiliary lemmas

**Lemma 9** Let $-\infty < \mu < \infty$, $0 < r \leq q \leq \infty$, and $b \in SV$. Then,

$$\left\|s^{-\mu-1/q} b(s) g(s)\right\|_{q,(0,t)} \leq C \left\|s^{-\mu-1/r} b(s) g(s)\right\|_{r,(0,t)} \tag{15}$$

for all $t > 0$ and all (Lebesgue-) measurable non-negative non-increasing functions $g$ on $(0,\infty)$, and

$$\left\|s^{-\mu-1/q} b(s) g(s)\right\|_{q,(t,\infty)} \leq C \left\|s^{-\mu-1/r} b(s) g(s)\right\|_{r,(t,\infty)}$$

for all $t > 0$ and all (Lebesgue-) measurable non-negative non-decreasing functions $g$ on $(0,\infty)$. The constants $C$ in both inequalities do not depend on $t$ and the function $g$.

*Proof.* Below in this proof, we denote by $C$ possibly different constants which do not depend on $t$ and $g$. We prove the first estimate. The second one can be proved similarly. Because the function $g$ is non-increasing, by Lemma 5 (ii) for all $s < t$, we have

$$s^{-\mu} b(s) g(s) \leq C g(s) \left\| u^{-\mu - 1/r} b(u) \right\|_{r,\left(\frac{s}{2},s\right)} \leq C \left\| u^{-\mu - 1/r} b(u) g(u) \right\|_{r,\left(\frac{s}{2},s\right)}$$

$$\leq C \left\| u^{-\mu - 1/r} b(u) g(u) \right\|_{r,(0,s)} \leq C \left\| u^{-\mu - 1/r} b(u) g(u) \right\|_{r,(0,t)}.$$

Thus, (15) is proved for $q=\infty$:

$$\operatorname*{esssup}_{0<s<t} s^{-\mu} b(s) g(s) \leq C \left\| s^{-\mu - 1/r} b(s) g(s) \right\|_{r,(0,t)}.$$

If $r < q < \infty$, using the last estimate, we get

$$\left\| s^{-\mu - 1/q} b(s) g(s) \right\|_{q,(0,t)} \leq \left( \operatorname*{esssup}_{0<s<t} s^{-\mu} b(s) g(s) \right)^{\frac{q-r}{q}} \left\{ \int_0^t \left( s^{-\mu} b(s) g(s) \right)^r \frac{ds}{s} \right\}^{\frac{1}{q}}$$

$$\leq C \left\| s^{-\mu - 1/r} b(s) g(s) \right\|_{r,(0,t)}.$$

This completes the proof. □

## 3. Main Results

Recall that the Fourier transform and its inverse may be uniquely extended to operators which act in $S'$. Note that $L_q \not\subset S'$ if $0<q<1$, $L_1 \subset S'$, and $L_{q,c;\mathbf{B}} \subset S'$ if $1 < q \leq \infty$. For $0 < R < \infty$, we denote

$$B_R := \{\xi \in \mathbb{R}^n : |\xi| \leq R\} \quad \text{and} \quad L_{q,c;\mathbf{B}}^R = \{f \in L_{q,c;\mathbf{B}} \cap S' : \operatorname{supp}\hat{f} \subset B_R\}.$$

From Nikol'skii's results [24, 25] (see also [23]), the next theorem follows.

**Theorem 10** Let $1 \leq q < p \leq \infty$. If $f \in L_q^R$, then $f \in L_p$ and

$$\|f\|_p \leq C R^{n\left(\frac{1}{q} - \frac{1}{p}\right)} \|f\|_q,$$

where $C$ depends only on $p$ and $q$.

Everywhere below $(p, b, \mathbf{A}), (q, c, \mathbf{B}) \in [1, \infty] \times [1, \infty] \times \mathbb{R}^2$, $q \leq p$, and (5) holds. For $F$ and $G$ being positive functions, we write $F \prec G$ if $F \leq CG$, where the constant $C$ is independent on all significant quantities such as $R$ and a function $f$. Such quantities as $p, b, \mathbf{A}, q, c, \mathbf{B}$ we consider as non-significant. Two positive functions $F$ and $G$ are considered equivalent ($F \approx G$) if $F \prec G$ and $G \prec F$. By $C$ we denote different constants which depend only on $p, b, \mathbf{A}, q, c$, and $\mathbf{B}$ and do not depend on $f$ and $R$. For shortness of the next formulations, let us denote

$$\mathbb{F} := \{(q, c, \mathbf{B}) : q = c = 1, \mathbf{B} = (0,0); \text{ or } 1 < q < \infty, 0 < c \leq \infty, \mathbf{B} = (\beta_0, \beta_\infty) \in \mathbb{R}^2\}.$$

Our goal is to obtain Nikol'skii–type inequalities for the functions $f \in L_{q,c;\mathbf{B}}^R$ of the form

$$\|f\|_{p,b;\mathbf{A}} \leq CG(R, p, b, \mathbf{A}, q, c, \mathbf{B}) \|f\|_{q,c;\mathbf{B}},$$

where $G(R, p, b, \mathbf{A}, q, c, \mathbf{B})$ is some expression depending on $R, p, b, \mathbf{A}, q, c, \mathbf{B}$. Or, writing it shorter,

$$\|f\|_{p,b;\mathbf{A}} \prec G(R, p, b, \mathbf{A}, q, c, \mathbf{B}) \|f\|_{q,c;\mathbf{B}}.$$

Obviously, only the situations $L_{q,c;\mathbf{B}} \not\subset L_{p,b;\mathbf{A}}$ are of interest, otherwise $G = 1$. Depending on $q$ and $p$, we consider three cases and use different techniques.

### 3.1. Case $1 < q < \infty$ and $q < p \leq \infty$

The next theorem extends Theorem 10 to the Lorentz–Zygmund spaces.

**Theorem 11** Let $(q, c, \mathbf{B}) \in \mathbb{F}$. For the triple $(p, b, \mathbf{A})$, we assume that either $q<p<\infty$ and $0 < b \leq \infty$ or $p=b=\infty$, $\mathbf{A} = 0$. Then for all $R > 0$ and $f \in L^R_{q,c;\mathbf{B}}$

$$\|f\|_{p,b;\mathbf{A}} \prec R^{n\left(\frac{1}{q}-\frac{1}{p}\right)} l^{(\overline{\mathbf{A}}-\overline{\mathbf{B}})}(R) \|f\|_{q,c;\mathbf{B}}.$$

***Proof.*** The proof is divided into four steps. In Step 1 we consider the case $1 < q < p < \infty$ and $b = \infty$. Here we use the technique from the proof of [18, Lemma 2.4]. In Step 2 we consider the case $1 < q < p < \infty$ and $0 < b < \infty$. Here we use some facts from the interpolation theory. In Step 3 and Step 4 we consider the last cases $(q, c, \mathbf{B}) = (1,1,0)$ and $(p, b, \mathbf{A}) = (\infty, \infty, 0)$, respectively, and use the results from above.

**Step 1.** First, consider the case $1 < q < p < \infty$ and $b = \infty$. Let $\chi \in C^\infty[0, \infty)$ is such that $\chi(u) = 1$ for $0 \leq u \leq 1$ and $\chi(u) = 0$ for $u \geq 2$. Define $v_R(x) := \mathcal{F}^{-1}[\chi(|\xi|/R)](x)$. Then for all $x \in \mathbb{R}^n$ and all $t > 0$

$$|v_R(x)| \prec \frac{R^n}{(1+R|x|)^n}, \quad v_R^*(t) \prec \frac{R^n}{(1+Rt^{1/n})^n}, \quad \text{and } v_R^{**}(t) \prec \min\left\{R^n, \frac{1}{t}\right\} \quad (16)$$

Because $\operatorname{supp} \hat{f} \subset B_R$ and $\operatorname{supp} \widehat{v_R} \subset B_{2R}$, we have $v_R * f = f$. Therefore, by O'Neil's inequality [26],

$$f^*(t) = (v_R * f)^*(t) \prec t v_R^{**}(t) f^{**}(t) + \int_t^\infty v_R^*(u) f^*(u) du. \quad (17)$$

By (17) and (16), it follows

$$\|f\|_{p,\infty;\mathbf{A}} = \sup_{0<t<\infty} \left( t^{\frac{1}{p}} l^{\mathbf{A}}(t) f^*(t) \right)$$

$$\prec \sup_{0<t<\infty} \left( t^{\frac{1}{p}+1} l^{\mathbf{A}}(t) v_R^{**}(t) f^{**}(t) \right) + \sup_{0<t<\infty} \left( t^{\frac{1}{p}} l^{\mathbf{A}}(t) \int_t^\infty v_R^*(u) f^*(u) du \right)$$

$$\prec \sup_{0<t<\infty} \left( t^{\frac{1}{p}+1} l^{\mathbf{A}}(t) \min\left\{R^n, \frac{1}{t}\right\} f^{**}(t) \right) + \sup_{0<t<\infty} \left( t^{\frac{1}{p}} l^{\mathbf{A}}(t) \int_t^\infty v_R^*(u) f^*(u) du \right)$$

$$:= S_1 + S_2.$$

Consider $S_1$.

$$S_1 \prec R^n \sup_{0<t<R^{-n}} \left( t^{\frac{1}{p}+1} l^{\mathbf{A}}(t) f^{**}(t) \right) + \sup_{R^{-n}<t<\infty} \left( t^{\frac{1}{p}} l^{\mathbf{A}}(t) f^{**}(t) \right) := S_{1,1} + S_{1,2}.$$

Since $1 < q < p$, we have $\frac{1}{p} - \frac{1}{q} + 1 > 0$ and $\frac{1}{p} - \frac{1}{q} < 0$. Thus, $t^{\frac{1}{p}-\frac{1}{q}+1} l^{\mathbf{A}-\mathbf{B}}(t)$ is equivalent to an increasing function and $t^{\frac{1}{p}-\frac{1}{q}} l^{\mathbf{A}-\mathbf{B}}(t)$ is equivalent to a decreasing function. Hence,

$$S_{1,1} = R^n \sup_{0<t<R^{-n}} \left( \left\{ t^{\frac{1}{p}-\frac{1}{q}+1} l^{\mathbf{A}-\mathbf{B}}(t) \right\} t^{\frac{1}{q}} l^{\mathbf{B}}(t) f^{**}(t) \right)$$

$$\prec R^{n\left(\frac{1}{q}-\frac{1}{p}\right)} l^{\mathbf{A}-\mathbf{B}}(R^{-n}) \sup_{0<t<R^{-n}} \left( t^{\frac{1}{q}} l^{\mathbf{B}}(t) f^{**}(t) \right)$$

$$\approx R^{n\left(\frac{1}{q}-\frac{1}{p}\right)} l^{(\overline{\mathbf{A}}-\overline{\mathbf{B}})}(R) \sup_{0<t<R^{-n}} \left( t^{\frac{1}{q}} l^{\mathbf{B}}(t) f^{**}(t) \right)$$

and

$$S_{1,2} = \sup_{R^{-n}<t<\infty} \left( \left\{ t^{\frac{1}{p}-\frac{1}{q}} l^{\mathbf{A}-\mathbf{B}}(t) \right\} t^{\frac{1}{q}} l^{\mathbf{B}}(t) f^{**}(t) \right)$$

$$\prec R^{n\left(\frac{1}{q}-\frac{1}{p}\right)} l^{(\overline{\mathbf{A}}-\overline{\mathbf{B}})}(R) \sup_{R^{-n}<t<\infty} \left( t^{\frac{1}{q}} l^{\mathbf{B}}(t) f^{**}(t) \right).$$

Therefore,

$$S_1 \prec R^{n\left(\frac{1}{q}-\frac{1}{p}\right)} l^{(\overline{\mathbf{A}}-\overline{\mathbf{B}})}(R) \left( \sup_{0<t<R^{-n}} \left( t^{\frac{1}{q}} l^{\mathbf{B}}(t) f^{**}(t) \right) + \sup_{R^{-n}<t<\infty} \left( t^{\frac{1}{q}} l^{\mathbf{B}}(t) f^{**}(t) \right) \right)$$

$$\approx R^{n\left(\frac{1}{q}-\frac{1}{p}\right)} l^{(\overline{\mathbf{A}}-\overline{\mathbf{B}})}(R) \sup_{0<t<\infty} \left( t^{\frac{1}{q}} l^{\mathbf{B}}(t) f^{**}(t) \right).$$

Since $f^{**}(t)$ is non-increasing, using Lemma 9 and [27, Theorem 3.8 (i)], we get

$$S_1 \prec R^{n\left(\frac{1}{q}-\frac{1}{p}\right)} l^{(\overline{\mathbf{A}}-\overline{\mathbf{B}})}(R) \left\| t^{\frac{1}{q}-\frac{1}{c}} l^{\mathbf{B}}(t) f^{**}(t) \right\|_{c,(0,\infty)} \approx R^{n\left(\frac{1}{q}-\frac{1}{p}\right)} l^{(\overline{\mathbf{A}}-\overline{\mathbf{B}})}(R) \|f\|_{q,c;\mathbf{B}}.$$

Consider $S_2$. Because $\frac{1}{p} > 0$ and $v_R^*(t) f^*(t)$ is non-increasing, applying the Hardy-type inequality [19, Lemma 2.7 (ii)], (15), and (16), we consider

$$S_2 \prec \sup_{0<t<\infty} \left( t^{1+\frac{1}{p}} l^{\mathbf{A}}(t) v_R^*(t) f^*(t) \right) \prec \left\| t^{1+\frac{1}{p}-\frac{1}{c}} l^{\mathbf{A}}(t) v_R^*(t) f^*(t) \right\|_{c,(0,\infty)}$$

$$\prec R^n \left\| t^{1+\frac{1}{p}-\frac{1}{c}} l^{\mathbf{A}}(t) \frac{1}{\left(1+Rt^{\frac{1}{n}}\right)^n} f^*(t) \right\|_{c,(0,\infty)}$$

$$\approx R^n \left( \left\| t^{1+\frac{1}{p}-\frac{1}{c}} l^{\mathbf{A}}(t) \frac{1}{\left(1+Rt^{\frac{1}{n}}\right)^n} f^*(t) \right\|_{c,(0,R^{-n})} + \left\| t^{1+\frac{1}{p}-\frac{1}{c}} l^{\mathbf{A}}(t) \frac{1}{\left(1+Rt^{\frac{1}{n}}\right)^n} f^*(t) \right\|_{c,(R^{-n},\infty)} \right).$$

Observing that

$$\left(1+Rt^{\frac{1}{n}}\right)^n \approx \begin{cases} 1, & 0<t<R^{-n}, \\ R^n t, & t \geq R^{-n}, \end{cases}$$

we conclude

$$S_2 \prec R^n \left\| t^{1+\frac{1}{p}-\frac{1}{c}} l^{\mathbf{A}}(t) f^*(t) \right\|_{c,(0,R^{-n})} + \left\| t^{\frac{1}{p}-\frac{1}{c}} l^{\mathbf{A}}(t) f^*(t) \right\|_{c,(R^{-n},\infty)} := S_{2,1} + S_{2,2}.$$

Since $t^{\frac{1}{p}-\frac{1}{q}+1} l^{\mathbf{A}-\mathbf{B}}(t)$ is equivalent to an increasing function and $t^{\frac{1}{p}-\frac{1}{q}} l^{\mathbf{A}-\mathbf{B}}(t)$ is equivalent to a decreasing function, we obtain

$$S_{2,1} = R^n \left\| \left\{ t^{1+\frac{1}{p}-\frac{1}{q}} l^{\mathbf{A}-\mathbf{B}}(t) \right\} t^{\frac{1}{q}-\frac{1}{c}} l^{\mathbf{B}}(t) f^*(t) \right\|_{c,(0,R^{-n})}$$

$$\prec R^{n\left(\frac{1}{q}-\frac{1}{p}\right)} l^{(\overline{\mathbf{A}}-\overline{\mathbf{B}})}(R) \left\| t^{\frac{1}{q}-\frac{1}{c}} l^{\mathbf{B}}(t) f^*(t) \right\|_{c,(0,R^{-n})}$$

and

$$S_{2,2} = \left\| t^{\frac{1}{p}-\frac{1}{c}} l^{\mathbf{A}}(t) f^*(t) \right\|_{c,(R^{-n},\infty)} = \left\| \left\{ t^{\frac{1}{p}-\frac{1}{q}} l^{\mathbf{A}-\mathbf{B}}(t) \right\} t^{\frac{1}{q}-\frac{1}{c}} l^{\mathbf{B}}(t) f^*(t) \right\|_{c,(R^{-n},\infty)}$$

$$\prec R^{n\left(\frac{1}{q}-\frac{1}{p}\right)} l^{(\overline{\mathbf{A}}-\overline{\mathbf{B}})}(R) \left\| t^{\frac{1}{q}-\frac{1}{c}} l^{\mathbf{B}}(t) f^*(t) \right\|_{c,(R^{-n},\infty)}.$$

Thus,

$$S_2 \approx R^{n\left(\frac{1}{q}-\frac{1}{p}\right)} l^{(\overline{\mathbf{A}}-\overline{\mathbf{B}})}(R) \left( \left\| t^{\frac{1}{q}-\frac{1}{c}} l^{\mathbf{B}}(t) f^*(t) \right\|_{c,(0,R^{-n})} + \left\| t^{\frac{1}{q}-\frac{1}{c}} l^{\mathbf{B}}(t) f^*(t) \right\|_{c,(R^{-n},\infty)} \right)$$

$$\approx R^{n\left(\frac{1}{q}-\frac{1}{p}\right)} l^{(\overline{\mathbf{A}}-\overline{\mathbf{B}})}(R) \left\| t^{\frac{1}{q}-\frac{1}{c}} l^{\mathbf{B}}(t) f^*(t) \right\|_{c,(0,\infty)} = R^{n\left(\frac{1}{q}-\frac{1}{p}\right)} l^{(\overline{\mathbf{A}}-\overline{\mathbf{B}})}(R) \|f\|_{q,c;\mathbf{B}}.$$

**Step 2.** Let now, $1 < q < p < \infty$ and $0 < b < \infty$. Take some $p_0$ such that $p < p_0 < \infty$. Because $q < p < p_0$, it is possible to choose $\theta$ ($0 < \theta < 1$) and $\mathbf{\Gamma}$ so that

$$\frac{1}{p} = \frac{1-\theta}{q} + \frac{\theta}{p_0} \quad \text{and} \quad \mathbf{A} = (1-\theta)\mathbf{B} + \theta\mathbf{\Gamma}.$$

This is equivalent to

$$\theta\left(\frac{1}{q}-\frac{1}{p_0}\right) = \frac{1}{q}-\frac{1}{p} \qquad \text{and} \qquad \theta(\mathbf{\Gamma}-\mathbf{B}) = \mathbf{A}-\mathbf{B}. \tag{18}$$

From Step 1, we have $f \in L^R_{p_0,\infty;\mathbf{\Gamma}}$ and

$$\|f\|_{p_0,\infty;\mathbf{\Gamma}} \prec R^{n\left(\frac{1}{q}-\frac{1}{p_0}\right)} l^{(\overline{\mathbf{\Gamma}}-\overline{\mathbf{B}})}(R) \|f\|_{q,c;\mathbf{B}}.$$

By (6), we have $L_{p,b;\mathbf{A}} = \left(L_{q,c;\mathbf{B}}, L_{p_0,\infty;\mathbf{\Gamma}}\right)_{\theta,b}$. Since $f \in L_{q,c;\mathbf{B}} \cap L_{p_0,\infty;\mathbf{\Gamma}}$, using now Lemma 1 (ii) and (18), we obtain

$$\|f\|_{p,b;\mathbf{A}} \prec \|f\|_{q,c;\mathbf{B}}^{(1-\theta)} \|f\|_{p_0,\infty;\mathbf{\Gamma}}^{\theta} \prec R^{\theta n\left(\frac{1}{q}-\frac{1}{p_0}\right)} l^{\theta(\overline{\mathbf{\Gamma}}-\overline{\mathbf{B}})}(R) \|f\|_{q,c;\mathbf{B}} = R^{n\left(\frac{1}{q}-\frac{1}{p}\right)} l^{(\overline{\mathbf{A}}-\overline{\mathbf{B}})}(R) \|f\|_{q,c;\mathbf{B}}.$$

**Step 3**. Let $q = c = 1, \mathbf{B} = 0, 1 < p < \infty$ and $0 < b \le \infty$. Take some $q_0$ such that $1 < q_0 < p$. By hypothesis $f \in L^R_1$. So, Theorem 10 implies that $f \in L^R_{q_0}$ and

$$\|f\|_{q_0} \prec R^{n\left(1-\frac{1}{q_0}\right)} \|f\|_1.$$

Using additionally the result either from Step 1 or from Step 2, we arrive at

$$\|f\|_{p,b;\mathbf{A}} \prec R^{n\left(\frac{1}{q_0}-\frac{1}{p}\right)} l^{\overline{\mathbf{A}}}(R) \|f\|_{q_0} \prec R^{n\left(1-\frac{1}{p}\right)} l^{\overline{\mathbf{A}}}(R) \|f\|_1.$$

**Step 4**. Let $p=b=\infty$ and $\mathbf{A} = 0$. Take some $p_0$ such that $q < p_0 < \infty$. From the previous steps, we have $f \in L^R_{p_0}$ and

$$\|f\|_{p_0} \prec R^{n\left(\frac{1}{q}-\frac{1}{p_0}\right)} l^{-\overline{\mathbf{B}}}(R) \|f\|_{q,c;\mathbf{B}}.$$

By Theorem 10, we conclude now that $f \in L_\infty$ and

$$\|f\|_\infty \prec R^{\frac{n}{p_0}} \|f\|_{p_0} \prec R^{\frac{n}{q}} l^{-\overline{\mathbf{B}}}(R) \|f\|_{q,c;\mathbf{B}}.$$

This completes the proof. $\qquad \square$

**Remark 12** Theorem 11 extends Lemma 2.4 in [18]. Compare also with Theorems 3 and 4 in [23], estimate (7) in [33, p. 31] applied with $d=0$ and $\alpha=0$, Theorem 6.5 in [8], and Theorem 5.2 in [6]. The limiting cases considered in the next two subsections are new, as far as we are aware.

### 3.2. Case $p = \infty$

In this subsection, we examine the case $p = \infty$, $0 < b \le \infty$, and $\mathbf{A} \ne \mathbf{0}$.

**Theorem 13** Let $(q,c,\mathbf{B}) \in \mathbb{F}$, $0 < b \le \infty$, and $\alpha_0 < -\frac{1}{b} < \alpha_\infty$. Then for all $R > 0$ and $f \in L^R_{q,c;\mathbf{B}}$

$$\|f\|_{\infty,b;\mathbf{A}} \prec R^{\frac{n}{q}} l^{\left(\overline{\mathbf{A}}-\overline{\mathbf{B}}+\frac{1}{b}\right)}(R) \|f\|_{q,c;\mathbf{B}}. \tag{19}$$

*Proof.* Note, that due to Theorem 11, if $f \in L^R_{q,c;\mathbf{B}}$, then $f \in L_\infty$. Hence $f \in L_{q,c;\mathbf{B}} \cap L_\infty$ and (7) implies $f \in L_{\infty,b;\mathbf{A}}$.

**Step 1**. Let $(q,c,\mathbf{B}) = (1,1,0)$ or $1 < q < 2$. First, we prove that if $\alpha_0 < -1 < \alpha_\infty$, then

$$\|f\|_{\infty,1;\mathbf{A}} \prec R^{\frac{n}{q}} l^{(\overline{\mathbf{A}}-\overline{\mathbf{B}}+1)}(R) \|f\|_{q,c;\mathbf{B}}. \tag{20}$$

If $1 < q < 2$, due to (13), Lemma 7, and (14), we have

$$f \in L^{B_R}_{q,c;\mathbf{B}} \xrightarrow{\mathcal{F}} L_{q',c;\overline{\mathbf{B}}}(B_R) \subset L_{1,1;\overline{\mathbf{A}}+1}(B_R) \xrightarrow{\mathcal{F}^{-1}} L_{\infty,1;\mathbf{A}},$$

where the norm of the embedding is bounded from above by $CR^{\frac{n}{q}} l^{(\overline{\mathbf{A}}-\overline{\mathbf{B}}+1)}(R)$. Similarly, in the case $(q,c,\mathbf{B}) = (1,1,0)$, by (12), it follows

$$f \in L_1^{B_R} \xrightarrow{\mathcal{F}} L_\infty(B_R) \subset L_{1,1;\overleftarrow{\mathbf{A}}+1}(B_R) \xrightarrow{\mathcal{F}^{-1}} L_{\infty,1;\mathbf{A}},$$

where the norm of the embedding is bounded from above by $CR^n l^{(\overleftarrow{\mathbf{A}}+1)}(R)$. Hence (20) holds. For arbitrary $b$ and $\mathbf{A}$, we put $\mathbf{\Gamma} = 2\left(\mathbf{A} + \frac{1}{b}\right) - 1$. According to [17, Theorem 7.9] and (7), we deduce that

$$L_{\infty,b;\mathbf{A}} = \left(L_{\infty,1;\mathbf{\Gamma}}, L_\infty\right)_{\frac{1}{2},b}.$$

Since $\gamma_0 < -1 < \gamma_\infty$, (20) implies that $f \in L_{\infty,1;\mathbf{\Gamma}}$. Hence $f \in L_{\infty,1;\mathbf{\Gamma}} \cap L_\infty$. Using now Lemma 1 (ii), (20), and Theorem 11, we obtain

$$\|f\|_{\infty,b;\mathbf{A}} \prec \|f\|_{\infty,1;\mathbf{\Gamma}}^{\frac{1}{2}} \|f\|_\infty^{\frac{1}{2}} \prec \left\{R^{\frac{n}{q}} l^{(\overleftarrow{\mathbf{\Gamma}}-\overleftarrow{\mathbf{B}}+1)}(R) R^{\frac{n}{q}} l^{-\overleftarrow{\mathbf{B}}}(R)\right\}^{\frac{1}{2}} \|f\|_{q,c;\mathbf{B}}$$
$$= R^{\frac{n}{q}} l^{\left(\overleftarrow{\mathbf{A}}-\overleftarrow{\mathbf{B}}+\frac{1}{b}\right)}(R) \|f\|_{q,c;\mathbf{B}}. \tag{21}$$

**Step 2**. In the case $2 < q < \infty$, we use one trick from [23]. Let $\rho(= 2, 3, \ldots)$ be such that $1 < \frac{q}{\rho} < 2$. Note that $\mathcal{F}\{f^\rho\} = \mathcal{F}\{f\} * \ldots * \mathcal{F}\{f\}$ and hence $\text{supp}(\mathcal{F}\{f^\rho\}) \subset B_{\rho R}$. Thus, by Lemma 8 and (21), we conclude that $f^\rho \in L_{\frac{q}{\rho},\frac{c}{\rho};\rho\mathbf{B}}^{\rho R}$ and

$$\|f\|_{\infty,b,\mathbf{A}} = \left\{\|f^\rho\|_{\infty,\frac{b}{\rho};\rho\mathbf{A}}\right\}^{\frac{1}{\rho}} \prec \left\{R^{\frac{n\rho}{q}} l^{\rho\left(\overleftarrow{\mathbf{A}}-\overleftarrow{\mathbf{B}}+\frac{1}{b}\right)}(R) \|f^\rho\|_{\frac{q}{\rho},\frac{c}{\rho},\rho\mathbf{B}}\right\}^{\frac{1}{\rho}} = R^{\frac{n}{q}} l^{\left(\overleftarrow{\mathbf{A}}-\overleftarrow{\mathbf{B}}+\frac{1}{b}\right)}(R) \|f\|_{q,c,\mathbf{B}}.$$

**Step 3**. It remains only to consider the case when $q = 2$. Take $q_0 = \frac{5}{2}$. From Theorem 11, we know that $f \in L_{q_0,c;\mathbf{B}}^R$ and

$$\|f\|_{q_0,c;\mathbf{B}} \prec R^{n\left(\frac{1}{q}-\frac{1}{q_0}\right)} \|f\|_{q,c,\mathbf{B}}.$$

From Step 2, we have

$$\|f\|_{\infty,b,\mathbf{A}} \prec R^{\frac{n}{q_0}} l^{\left(\overleftarrow{\mathbf{A}}-\overleftarrow{\mathbf{B}}+\frac{1}{b}\right)}(R) \|f\|_{q_0,c,\mathbf{B}}.$$

Combining the two last inequalities, we get (19) in this case. □

### 3.3. Case $q = p$

The next two theorems examine the limiting case $q = p$.

**Theorem 14** Let $(q, c, \mathbf{B}) \in \mathbb{F}$, $0 < b \leq c \leq \infty$, and $\alpha_\infty + \frac{1}{b} < \beta_\infty + \frac{1}{c}$. Put $\mathbf{\Gamma} = \overleftarrow{\mathbf{A}} + \frac{1}{b} - \overleftarrow{\mathbf{B}} - \frac{1}{c}$. In addition, assume that one of the following conditions is met

(i) $\alpha_0 + \frac{1}{b} > \beta_0 + \frac{1}{c}$ and $\mathbf{\Delta} = 0$;

(ii) $\alpha_0 + \frac{1}{b} = \beta_0 + \frac{1}{c}$ and $\mathbf{\Delta} = \left(0, \frac{1}{b} - \frac{1}{c}\right)$.

Then for all $R > 0$ and $f \in L_{q,c;\mathbf{B}}^R$

$$\|f\|_{q,b;\mathbf{A}} \prec l^{\mathbf{\Gamma}}(R) ll^{\mathbf{\Delta}}(R) \|f\|_{q,c;\mathbf{B}}. \tag{22}$$

*Proof*. (Cf. [11, proof of Theorem 3].) We examine only the case $0 < b < c < \infty$. For all other cases, the proofs are similar. The assertion $f \in L_{q,b;\mathbf{A}}$ follows from the reasonings below. We write

$$\|f\|_{q,b;\mathbf{A}}^b = \int_0^\infty \left(t^{\frac{1}{q}} l^{\mathbf{A}}(t) f^*(t)\right)^b \frac{dt}{t} = \int_0^{R^{-n}} \ldots + \int_{R^{-n}}^\infty \ldots := I_1(R) + I_2(R).$$

First, we estimate $I_1$:

$$I_1(R) = \int_0^{R^{-n}} t^{\frac{b}{2q}} \left( t^{\frac{1}{2q}} l^{\mathbf{A}}(t) f^*(t) \right)^b \frac{dt}{t} \leq R^{-\frac{bn}{2q}} \|f\|_{2q,b;\mathbf{A}}^b.$$

Using Theorem 11, we get $f \in L_{2q,c;\mathbf{B}}^R$ and

$$I_1(R) \prec R^{-\frac{bn}{2q}} \left\{ R^{n\left(\frac{1}{q}-\frac{1}{2q}\right)} l^{(\overleftarrow{\mathbf{A}}-\overleftarrow{\mathbf{B}})}(R) \|f\|_{q,c;\mathbf{B}} \right\}^b = \left\{ l^{(\overleftarrow{\mathbf{A}}-\overleftarrow{\mathbf{B}})}(R) \|f\|_{q,c;\mathbf{B}} \right\}^b.$$

Note that $b < c$ implies $\overleftarrow{\mathbf{A}} - \overleftarrow{\mathbf{B}} < \boldsymbol{\Gamma}$ and $l^{(\overleftarrow{\mathbf{A}}-\overleftarrow{\mathbf{B}})}(R) < l^{\boldsymbol{\Gamma}}(R)$. Therefore, to prove each of the statements of Theorem 14, it suffices to show that

$$I_2(R) \prec \left\{ l^{\boldsymbol{\Gamma}}(R) l l^{\boldsymbol{\Delta}}(R) \|f\|_{q,c;\mathbf{B}} \right\}^b.$$

Using Hölder's inequality with exponents $\frac{c}{b}$ and $\frac{c}{c-b}$, we derive that

$$I_2(R) \leq \left( \int_{R^{-n}}^\infty \left( l^{\mathbf{A}-\mathbf{B}}(t) \right)^{\frac{bc}{c-b}} \frac{dt}{t} \right)^{\frac{c-b}{c}} \left( \int_{R^{-n}}^\infty \left( t^{\frac{1}{q}} l^{\mathbf{B}}(t) f^*(t) \right)^c \frac{dt}{t} \right)^{\frac{b}{c}}$$

$$\leq \left( \int_{R^{-n}}^\infty \left( l^{\mathbf{A}-\mathbf{B}}(t) \right)^{\frac{bc}{c-b}} \frac{dt}{t} \right)^{\frac{c-b}{c}} \|f\|_{q,c;\mathbf{B}}^b.$$

The condition $\alpha_\infty + \frac{1}{b} < \beta_\infty + \frac{1}{c}$ is equivalent to $(\alpha_\infty - \beta_\infty) \frac{bc}{c-b} < -1$ and the integral above converges. Hence, for all $R \leq 1$

$$\int_{R^{-n}}^\infty \left( l^{\mathbf{A}-\mathbf{B}}(t) \right)^{\frac{bc}{c-b}} \frac{dt}{t} = \int_{R^{-n}}^\infty \left( l^{\alpha_\infty-\beta_\infty}(t) \right)^{\frac{bc}{c-b}} \frac{dt}{t} \approx l^{\frac{bc}{c-b}\left(\alpha_\infty+\frac{1}{b}-\beta_\infty-\frac{1}{c}\right)}(R)$$

and

$$I_2(R) \prec l^{b\left(\alpha_\infty+\frac{1}{b}-\beta_\infty-\frac{1}{c}\right)}(R) \|f\|_{q,c;\mathbf{B}}^b.$$

Thus, the estimate (22) holds for $R \leq 1$. Let now $R > 1$. Because $l^{\boldsymbol{\Gamma}}(R) l l^{\boldsymbol{\Delta}}(R) \to \infty$ if $R \to \infty$, it is enough to suppose that $R$ is large enough, e.g. $R > 2$. If the condition (i) holds then $(\alpha_0 - \beta_0) \frac{bc}{c-b} > -1$ and

$$\int_{R^{-n}}^\infty \left( l^{\mathbf{A}-\mathbf{B}}(t) \right)^{\frac{bc}{c-b}} \frac{dt}{t} \approx \int_{R^{-n}}^1 \left( l^{\alpha_0-\beta_0}(t) \right)^{\frac{bc}{c-b}} \frac{dt}{t}$$

$$\approx l^{(\alpha_0-\beta_0)\frac{bc}{c-b}+1}(R) = l^{\frac{bc}{c-b}\left(\alpha_0+\frac{1}{b}-\beta_0-\frac{1}{c}\right)}(R).$$

So, in the case $\alpha_0 + \frac{1}{b} > \beta_0 + \frac{1}{c}$,

$$I_2(R) \prec l^{b\left(\alpha_0+\frac{1}{b}-\beta_0-\frac{1}{c}\right)}(R) \|f\|_{q,c;\mathbf{B}}^b.$$

If the condition (ii) holds then $(\alpha_0 - \beta_0) \frac{bc}{c-b} = -1$. Hence,

$$\int_{R^{-n}}^\infty \left( l^{\mathbf{A}-\mathbf{B}}(t) \right)^{\frac{bc}{c-b}} \frac{dt}{t} \approx \int_{R^{-n}}^1 l^{-1}(t) \frac{dt}{t} \approx ll(R)$$

and

$$I_2(R) \prec l l^{b\left(\frac{1}{b}-\frac{1}{c}\right)}(R) \|f\|_{q,c;\mathbf{B}}^b.$$

This completes the proof. □

**Theorem 15** Let $(q, c, \mathbf{B}) \in \mathbb{F}$, $0 < c \leq b \leq \infty$, $\alpha_0 \geq \beta_0$, and $\alpha_\infty < \beta_\infty$. Then for all $R > 0$ and $f \in L_{q,c;\mathbf{B}}^R$

$$\|f\|_{q,b;\mathbf{A}} \prec l^{(\overleftarrow{\mathbf{A}}-\overleftarrow{\mathbf{B}})}(R) \|f\|_{q,c;\mathbf{B}}.$$

*Proof.* Using Theorem 14 for the case $c = b$, we get

$$\|f\|_{q,b;\mathbf{A}} \prec l^{(\overleftarrow{\mathbf{A}}-\overleftarrow{\mathbf{B}})}(R) \|f\|_{q,b;\mathbf{B}}.$$

If $c < b$, then [27, Theorem 4.1] $L_{q,c;\mathbf{B}} \subset L_{q,b;\mathbf{B}}$. Consequently

$$\|f\|_{q,b;\mathbf{A}} \prec l^{(\overline{\mathbf{A}}-\overline{\mathbf{B}})}(R)\|f\|_{q,b;\mathbf{B}} \prec l^{(\overline{\mathbf{A}}-\overline{\mathbf{B}})}(R)\|f\|_{q,c;\mathbf{B}}.$$

This completes the proof. □

## 4. Applications

In this section we define some Besov–type spaces (as well homogeneous counterparts) of logarithmic smoothness based on Lorentz–Zygmund spaces and apply main results to establish corresponding embeddings between these spaces. Besov spaces play an important role in the theory of differentiable functions of several variables, approximation theory, the theory of partial differential equations, and other areas of mathematics. Many authors in different contexts have investigated Besov spaces of logarithmic smoothness. See, e.g. [3, 5, 12, 13, 14, 15, 18, 20, 21, 22, 30, 33, 35] and the references given therein.

Let $\varphi \in S(\mathbb{R}^n)$ be such function that [3, Lemma 6.1.7]
$$\begin{cases} \operatorname{supp} \varphi = \{\xi \colon 2^{-1} \leq |\xi| \leq 2\}, \\ \varphi(\xi) > 0 \text{ for } 2^{-1} < |\xi| < 2, \\ \sum_{k=-\infty}^{\infty} \varphi(2^{-k}\xi) = 1, (\xi \neq 0). \end{cases}$$

Using function $\varphi$, we define functions $\varphi_k$ and $\psi$ by
$$\begin{cases} \varphi_k(\xi) = \varphi(2^{-k}\xi) \ (k = 0, \pm 1, \pm 2, \dots), \\ \psi(\xi) = 1 - \sum_{k=1}^{\infty} \varphi(2^{-k}\xi). \end{cases}$$

We define Besov–type spaces, which have classical smoothness $\sigma$ and additional logarithmic smoothness $\gamma$ (or $\mathbf{\Gamma}$ for homogeneous spaces).

**Definition 16** Let $-\infty<\sigma,\gamma<\infty$, $0<p,b,u\leq\infty$, $\mathbf{A} = (\alpha_0, \alpha_\infty)$, and the space $L_{p,b;\mathbf{A}}$ is not trivial. The Besov–type space $B^{\sigma,\gamma}_{(p,b;\mathbf{A}),u}$ based on $L_{p,b;\mathbf{A}}$ consists of all $f \in S'$ having a finite quasi-norm
$$\|f\|_{B^{\sigma,\gamma}_{(p,b;\mathbf{A}),u}} = \|\mathcal{F}^{-1}(\psi \mathcal{F}f)\|_{p,b;\mathbf{A}} + \left\{\sum_{k=1}^{\infty}\left(2^{\sigma k}(1+k)^{\gamma}\|\mathcal{F}^{-1}(\varphi_k \mathcal{F}f)\|_{p,b;\mathbf{A}}\right)^u\right\}^{\frac{1}{u}}$$
with the usual modification for $u=\infty$.

**Definition 17** Let $-\infty<\sigma<\infty$, $0<p,b,u\leq\infty$, $\mathbf{A} = (\alpha_0, \alpha_\infty)$, $\mathbf{\Gamma} = (\gamma_0, \gamma_\infty)$, and the space $L_{p,b;\mathbf{A}}$ is not trivial. The homogeneous Besov–type space $\dot{B}^{\sigma,\mathbf{\Gamma}}_{(p,b;\mathbf{A}),u}$ consists of all $f \in S'$ having a finite semi-(quasi-)norm
$$\|f\|_{\dot{B}^{\sigma,\mathbf{\Gamma}}_{(p,b;\mathbf{A}),u}} = \left\{\sum_{k=-\infty}^{\infty}\left(2^{\sigma k}l^{\mathbf{\Gamma}}(2^k)\|\mathcal{F}^{-1}(\varphi_k \mathcal{F}f)\|_{p,b;\mathbf{A}}\right)^u\right\}^{\frac{1}{u}}.$$
with the usual modification for $u=\infty$.

We denote the spaces $B^{\sigma,\gamma}_{(p,p;0),u}$ ($\dot{B}^{\sigma,\mathbf{\Gamma}}_{(p,p;0),u}$) by $B^{\sigma,\gamma}_{p,u}$ ($\dot{B}^{\sigma,\mathbf{\Gamma}}_{p,u}$). The spaces $B^{\sigma,0}_{p,u}$ coincide with the classical Besov spaces $B^{\sigma}_{p,u}$ (see, e.g. [3, 35]), and the spaces $\dot{B}^{\sigma,0}_{p,u}$ coincide with the classical homogeneous Besov spaces $\dot{B}^{\sigma}_{p,u}$ (see, e.g. [3, 6.3], [30, Definition 4.4], [20]).

Recall, $B_R$ is the ball in $\mathbb{R}^n$ of the radius $R$ and
$$\mathbb{F} := \{(q, c, \mathbf{B}) \colon q = c = 1, \mathbf{B} = 0; or\ 1 < q < \infty\}.$$
Let $f \in S'$ and $g := \mathcal{F}^{-1}(\psi \mathcal{F}f)$. Obviously, for the Fourier transform of the function $g$ it holds $supp \mathcal{F}g = supp(\psi \mathcal{F}f) \subset supp\,\psi = B_2$. Similarly, if $g_k := \mathcal{F}^{-1}(\varphi_k \mathcal{F}f)$ then $supp \mathcal{F}g_k = supp(\varphi_k \mathcal{F}f) \subset supp\,\varphi_k \subset B_{2^{k+1}}$. If additionally $(q, c, \mathbf{B}) \in \mathbb{F}$ and $g, g_k \in L_{q,c,\mathbf{B}}$ then for these functions we can use the results from Section 3. In this way, we get following corollaries.

**Corollary 18** (Cf. [18, Corollary 2.8].) Let $(q, c, \mathbf{B}) \in \mathbb{F}$, $-\infty<\sigma,\gamma<\infty$, and $0<u\leq\infty$. For the triple $(p, b, \mathbf{A})$, we assume that either $q<p<\infty$ or $p=b=\infty$, $\mathbf{A} = 0$. Then

$$B_{(q,c;\mathbf{B}),u}^{\sigma+n\left(\frac{1}{q}-\frac{1}{p}\right),\gamma+(\alpha_0-\beta_0)} \subset B_{(p,b;\mathbf{A}),u}^{\sigma,\gamma}.$$

In particular, (cf. [22, Proposition 1.9 (iv)])

$$B_{q,u}^{\sigma+n\left(\frac{1}{q}-\frac{1}{p}\right),\gamma} \subset B_{p,u}^{\sigma,\gamma}.$$

*Proof.* The Fourier transform of each of the function $\mathcal{F}^{-1}(\varphi_k \mathcal{F} f)$ ($k = 1, 2, \ldots$) has support which lies in $B_{2^{k+1}}$. Thus, $\mathcal{F}^{-1}(\varphi_k \mathcal{F} f) \in L_{q,c;\mathbf{B}}^{2^{k+1}}$ and by Theorem 11, we get

$$\|\mathcal{F}^{-1}(\varphi_k \mathcal{F} f)\|_{p,b;\mathbf{A}} \leq C 2^{kn\left(\frac{1}{q}-\frac{1}{p}\right)} (1+k)^{(\alpha_0-\beta_0)} \|\mathcal{F}^{-1}(\varphi_k \mathcal{F} f)\|_{q,c;\mathbf{B}}.$$

where the constant $C$ does not depend on $k$ and $f$. The analogous statement holds for the function $\mathcal{F}^{-1}(\psi \mathcal{F} f)$, too. Let $u<\infty$. (The case $u=\infty$ can be considered similarly.) Then

$$\|f\|_{B_{(p,b;\mathbf{A}),u}^{\sigma,\gamma}} = \|\mathcal{F}^{-1}(\psi \mathcal{F} f)\|_{p,b;\mathbf{A}} + \left\{\sum_{k=1}^{\infty} \left(2^{\sigma k}(1+k)^{\gamma} \|\mathcal{F}^{-1}(\varphi_k \mathcal{F} f)\|_{p,b;\mathbf{A}}\right)^u\right\}^{\frac{1}{u}}$$

$$\leq C\left(\|\mathcal{F}^{-1}(\psi \mathcal{F} f)\|_{q,c;\mathbf{B}} + \left\{\sum_{k=1}^{\infty} \left(2^{\left(\sigma+n\left(\frac{1}{q}-\frac{1}{p}\right)\right)k}(k+1)^{\gamma+(\alpha_0-\beta_0)}\|\mathcal{F}^{-1}(\varphi_k \mathcal{F} f)\|_{q,c;\mathbf{B}}\right)^u\right\}^{\frac{1}{u}}\right)$$

$$= C\|f\|_{B_{(q,c;\mathbf{B}),u}^{\sigma+n\left(\frac{1}{q}-\frac{1}{p}\right),\gamma+(\alpha_0-\beta_0)}}.$$

This completes the proof. □

The next corollaries can be proved analogously.

**Corollary 19** Let $(q, c, \mathbf{B}) \in \mathbb{F}$, $-\infty<\sigma<\infty$, and $0<u\leq\infty$. For the triple $(p, b, \mathbf{A})$, we assume that either $q<p<\infty$ or $p=b=\infty$, $\mathbf{A} = 0$. Then

$$\dot{B}_{(q,c;\mathbf{B}),u}^{\sigma+n\left(\frac{1}{q}-\frac{1}{p}\right),\mathbf{\Gamma}+(\overleftarrow{\mathbf{A}}-\overleftarrow{\mathbf{B}})} \subset \dot{B}_{(p,b;\mathbf{A}),u}^{\sigma,\mathbf{\Gamma}}.$$

In particular, (cf. [20, Theorem 2.1])

$$\dot{B}_{q,u}^{\sigma+n\left(\frac{1}{q}-\frac{1}{p}\right),\mathbf{\Gamma}} \subset \dot{B}_{p,u}^{\sigma,\mathbf{\Gamma}}.$$

**Corollary 20** Let $(q, c, \mathbf{B}) \in \mathbb{F}$, $-\infty<\sigma,\gamma<\infty$, $0<u\leq\infty$, $0<b\leq\infty$, and $\alpha_0 < -\frac{1}{b} < \alpha_\infty$. Then

$$B_{(q,c;\mathbf{B}),u}^{\sigma+\frac{n}{q},\gamma+\left(\alpha_0+\frac{1}{b}-\beta_0\right)} \subset B_{(\infty,b;\mathbf{A}),u}^{\sigma,\gamma}.$$

**Corollary 21** Let $(q, c, \mathbf{B}) \in \mathbb{F}$, $-\infty<\sigma<\infty$, $0<u\leq\infty$, $0<b\leq\infty$, and $\alpha_0 < -\frac{1}{b} < \alpha_\infty$. Then

$$\dot{B}_{(q,c;\mathbf{B}),u}^{\sigma+\frac{n}{q},\mathbf{\Gamma}+\left(\overleftarrow{\mathbf{A}}-\overleftarrow{\mathbf{B}}+\frac{1}{b}\right)} \subset \dot{B}_{(\infty,b;\mathbf{A}),u}^{\sigma,\mathbf{\Gamma}}.$$

**Corollary 22** Let $(q, c, \mathbf{B}) \in \mathbb{F}$, $-\infty<\sigma,\gamma<\infty$, $0<u\leq\infty$, $0<b\leq c\leq\infty$, $\alpha_\infty + \frac{1}{b} < \beta_\infty + \frac{1}{c}$, and $\alpha_0 + \frac{1}{b} > \beta_0 + \frac{1}{c}$. Then

$$B_{(q,c;\mathbf{B}),u}^{\sigma,\gamma+\left(\alpha_0+\frac{1}{b}-\beta_0-\frac{1}{c}\right)} \subset B_{(q,b;\mathbf{A}),u}^{\sigma,\gamma}.$$

**Corollary 23** Let $(q, c, \mathbf{B}) \in \mathbb{F}$, $-\infty<\sigma<\infty$, $0<u\leq\infty$, $0<b\leq c\leq\infty$, $\alpha_\infty + \frac{1}{b} < \beta_\infty + \frac{1}{c}$, and $\alpha_0 + \frac{1}{b} > \beta_0 + \frac{1}{c}$. Then

$$\dot{B}^{\sigma,\mathbf{\Gamma}+\overline{\mathbf{A}}-\overline{\mathbf{B}}+\frac{1}{b}-\frac{1}{c}}_{(q,c;\mathbf{B}),u} \subset \dot{B}^{\sigma,\mathbf{\Gamma}}_{(q,b;\mathbf{A}),u}.$$

**Corollary 24** Let $(q,c,\mathbf{B}) \in \mathbb{F}$, $-\infty<\sigma,\gamma<\infty$, $0<u\leq\infty$, $0<b\leq c\leq\infty$, $\alpha_\infty < \beta_\infty$, and $\alpha_0 \geq \beta_0$. Then
$$B^{\sigma,\gamma+(\alpha_0-\beta_0)}_{(q,c;\mathbf{B}),u} \subset B^{\sigma,\gamma}_{(q,b;\mathbf{A}),u}.$$

**Corollary 25** Let $(q,c,\mathbf{B}) \in \mathbb{F}$, $-\infty<\sigma<\infty$, $0<u\leq\infty$, $0<b\leq c\leq\infty$, $\alpha_\infty < \beta_\infty$, and $\alpha_0 \geq \beta_0$. Then
$$\dot{B}^{\sigma,\mathbf{\Gamma}+\overline{\mathbf{A}}-\overline{\mathbf{B}}}_{(q,c;\mathbf{B}),u} \subset \dot{B}^{\sigma,\mathbf{\Gamma}}_{(q,b;\mathbf{A}),u}.$$